\input amstex
\input amsppt.sty
\magnification=\magstep1
\hsize=33truecc
\vsize=22.2truecm
\baselineskip=16truept
\NoBlackBoxes
\nologo
\pageno=1
\topmatter
\TagsOnRight

\def\N{\Bbb N}
\def\Z{\Bbb Z}

\def\l{\left}
\def\r{\right}
\def\b{\bigg}

\def\({\b(}
\def\[{\b[}
\def\){\b)}
\def\]{\b]}

\def\t{\text}
\def\f{\frac}
\def\mo{\roman{mod}}

\def\em{\emptyset}
\def\se {\subseteq}

\def\sm{\setminus}

\def\eq{\equiv}

\def\ls{\leqslant}
\def\gs{\geqslant}

\def\da{\delta}

\def\Proof{\noindent{\it Proof}}
\def\Remark{\noindent{\it Remark}}

\def\Ack{\noindent {\bf Acknowledgments}}
\hbox{J. Number Theory 175(2017), 167--190.}
\medskip
\title Refining Lagrange's four-square theorem\endtitle
\author Zhi-Wei Sun \endauthor
\affil Department of Mathematics, Nanjing University
     \\Nanjing 210093, People's Republic of China
    \\  zwsun\@nju.edu.cn
    \\ {\tt http://math.nju.edu.cn/$\sim$zwsun}
 \endaffil
\abstract Lagrange's four-square theorem asserts that any $n\in\N=\{0,1,2,\ldots\}$ can be written as
the sum of four squares. This can be further refined in various ways. We show that
any $n\in\N$ can be written as $x^2+y^2+z^2+w^2$ with $x,y,z,w\in\Z$ such that $x+y+z$ (or $x+2y$, or $x+y+2z$) is a square (or a cube).
We also prove that any $n\in\N$ can be written as
$x^2+y^2+z^2+w^2$ with $x,y,z,w\in\N$ such that $P(x,y,z)$ is a square, whenever $P(x,y,z)$ is among the polynomials
$$\gather x,\ 2x,\ x-y,\ 2x-2y,\ a(x^2-y^2)\ (a=1,2,3),\ x^2-3y^2,\ 3x^2-2y^2,
\\ x^2+ky^2\, (k=2,3,5,6,8,12),\ (x+4y+4z)^2+(9x+3y+3z)^2,
\\ x^2y^2+y^2z^2+z^2x^2,\ x^4+8y^3z+8yz^3,\ x^4+16y^3z+64yz^3.
\endgather$$
We also pose some conjectures for further research; for example, our 1-3-5-Conjecture states that any $n\in\N$
can be written as $x^2+y^2+z^2+w^2$ with $x,y,z,w\in\N$ such that $x+3y+5z$ is a square.
\endabstract
\thanks 2010 {\it Mathematics Subject Classification}.
Primary 11E25; Secondary 11B75, 11D85, 11E20.
\newline \indent Supported by the National Natural Science
Foundation (Grant No. 11571162) of China.
\endthanks
\endtopmatter
\document

\heading{1. Introduction}\endheading

Let $\N=\{0,1,2,\ldots\}$ be the set of all natural numbers (nonnegative integers).
Lagrange's four-square theorem (cf. [N96, pp.\,5-7]) states that any $n\in\N$ can be written as $x^2+y^2+z^2+w^2$ with $x,y,z,w\in\N$.

It is known that for any $a,b,c\in\Z^+=\{1,2,3,\ldots\}$ the set
$$E(a,b,c):=\{n\in\N:\ n\not= ax^2+by^2+cz^2\ \t{for any}\ x,y,z\in\Z\}\tag1.1$$
is not only nonempty but also infinite. A classical theorem of Gauss and Legendre (cf. [N96, p.\,23] or [MW, p.\,42]) asserts that
$$E(1,1,1)=\{4^k(8l+7):\ k,l\in\N\}.\tag1.2$$

In this paper we study various refinements of Lagrange's theorem.

\proclaim{Theorem 1.1} Let $a\in\{1,4\}$ and $m\in\{4,5,6\}$. Then, any $n\in\N$ can be written as $ax^m+y^2+z^2+w^2$ with $x,y,z,w\in\N$.
\endproclaim
\Remark\ 1.1. See [S16, A270969, A273915 and A273429] for related data; for example,
$$71 = 1^4 + 3^2 + 5^2 + 6^2,\ 240=2^5 + 0^2 + 8^2 + 12^2\  \t{and}\ \ 624 = 2^6 + 4^2 + 12^2 + 20^2.$$
In addition, we conjecture that any $n\in\N$ can be written as $x^2+y^3+z^4+2w^4$ with $x,y,z,w\in\N$ (cf. [S16, A262827])
and that each $n\in\N$ can be written as $x^5+y^4+z^2+3w^2$ with $x,y,z,w\in\N$ (cf. [S16, A273917]).
\medskip

For convenience we introduce the following definition.
\medskip

{\it Definition} 1.1. A polynomial $P(x,y,z,w)$ with integer coefficients is called a {\it suitable} polynomial if any $n\in\N$
can be written as $x^2+y^2+z^2+w^2$ with $x,y,z,w\in\N$ such that $P(x,y,z,w)$ is a square.
\medskip

Theorem 1.1 with $a=1,4$ and $m=4$ indicates that both $x$ and $2x$ are suitable.
For any squarefree integer $a\gs 3$, the polynomial $ax$ is not suitable since $7\not=(ax^2)^2+y^2+z^2+w^2$
for all $x,y,z,w\in\N$.

\proclaim{Theorem 1.2} {\rm (i)} The polynomials $x-y$ and $2(x-y)$ are both suitable.

{\rm (ii)} Let $c\in\{1,2,4\}$. Then, any $n\in\N$ can be written as $x^2+y^2+z^2+w^2$  $(x,y,z,w\in\Z)$
with $x+y=ct^3$ for some $t\in\Z$.

{\rm (iii)} Let $d\in\{1,2\}$ and $m\in\{2,3\}$.
Then, any $n\in\N$ can be written as $x^2+y^2+z^2+w^2$ with $x,y,z,w\in\Z$ such that $x+2y=dt^m$ for some $t\in\Z$.

{\rm (iv)} Any $n\in\N$ can be written as $x^2+y^2+z^2+w^2$ with $x,y,z,w\in\Z$ such that $xy+yz+zw+wx=(x+z)(y+w)$ is a square.

{\rm (v)} Each $n\in\Z^+$ can be written as $4^k(1+4x^2+y^2)+z^2$ with $k,x,y,z\in\N$.
\endproclaim
\Remark\ 1.2. We even conjecture that any $n\in\N$ not of the form $2^{6k+3}\times7$ ($k\in\N$)
can be written as $x^2+y^2+z^2+w^2$ ($x,y,z,w\in\N$) with $x-y$ a cube, and that any $n\in\Z^+$
can be written as $4^k(1+4x^2+y^2)+z^2$ ($k,x,y,z\in\N$) with $x\ls y$ (or $x\ls z$).

\proclaim{Theorem 1.3} {\rm (i)} Let $m\in\{2,3\}$ and $(c,d)\in S(m)$, where
$$S(3)=\{(1,1),(1,2),(2,1),(2,2),(2,4)\}\ \ \t{and}\ \ S(2)=S(3)\cup\{(2,3),(2,6)\}.$$
Then any $n\in\N$ can be written as $x^2+y^2+z^2+w^2$ with $x,y,z,w\in\Z$ such that $x+y+cz=dt^m$
for some $t\in\Z$.

{\rm (ii)} If any odd integer $n>2719$ can be represented by $x^2 + y^2 + 10z^2$ with $x,y,z\in\Z$, then
each $n\in\N$ can be written as $x^2+y^2+z^2+w^2\ (x,y,z,w\in\Z)$ with $x+3y$ a square.

{\rm (iii)} If any integer $n\gs1190$ not divisible by $16$ can be written as $x^2+10y^2+(2z^2+5^3r^4)/7$ with $r\in\{0,1,2,3\}$ and $x,y,z\in\Z$, then
any $n\in\N$ can be written as $x^2+y^2+z^2+w^2$ $(x,y,z,w\in\Z)$ with $x+3y+5z$ a square.
\endproclaim
\Remark\ 1.3. Concerning the condition of Theorem 1.3(ii),
in 1916 S. Ramanujan [R] conjectured that any odd integer $n>2719$ can be represented by $x^2 + y^2 + 10z^2$ with $x,y,z\in\Z$,
and this was proved by K. Ono and K. Soundararajan [OS] in 1997 under the GRH (Generalized Riemann Hypothesis).
As for part (iii) of Theorem 1.3, we guess that any integer $n\gs1190$ not divisible by $16$ can be written as $x^2+10y^2+(2z^2+5^3r^4)/7$ with $r\in\{0,1,2,3\}$ and $x,y,z\in\Z$.
\medskip

\proclaim{Theorem 1.4} {\rm (i)} Any $n\in\N$ can be written as $x^2+y^2+z^2+w^2$ with $x,y,z,w\in\N$ such that
$P(x,y,z)=0$, whenever $P(x,y,z)$ is among the polynomials
$$\aligned &x(x-y),\ x(x-2y),\ (x-y)(x-2y),\ (x-y)(x-3y),
\\&x(x+y-z),\ (x-y)(x+y-z),\ (x-2y)(x+y-z).
\endaligned\tag1.3$$

{\rm (ii)} Any $n\in\N$ can be written as $x^2+y^2+z^2+w^2$ with $x,y,z,w\in\N$ such that
$(2x-3y)(x+y-z)=0$, provided that
$$\{x^2+y^2+13z^2:\ x,y,z\in\N\}\supseteq\{8q+5:\ q\in\N\}.\tag1.4$$

{\rm (iii)} Any $n\in\Z^+$ can be written as $x^2+y^2+z^2+w^2$ with $x,y,z,w\in\N$ and $z>0$ such that $(x-y)z$ is a square.

{\rm (iv)} Any $n\in\Z^+$ can be written as $x^2+y^2+z^2+w^2$ with $x,y,z,w\in\N$ and $y>0$ such that $x+4y+4z$ and $9x+3y+3z$
are the two legs of a right triangle with positive integer sides.

{\rm (v)} Each positive integer can be written as $x^2+y^2+z^2+w^2$ with $x,y,z\in\N$ and $w\in\Z^+$ such that $x^2y^2+y^2z^2+z^2x^2$ is a square.
Also, any $n\in\Z^+$ can be written as $x^2+y^2+z^2+w^2$ with $x,y,z\in\Z$ and $w\in\Z^+$ such that $x^2y^2+4y^2z^2+4z^2x^2$ is a square.

{\rm (vi)} Any $n\in\Z^+$ can be written as $x^2+y^2+z^2+w^2$ with $x,y,z\in\N$ and $w\in\Z^+$ such that $x^4+8y^3z+8yz^3$ is a fourth power.
Also, any $n\in\Z^+$ can be written as $x^2+y^2+z^2+w^2$ with $x,y,z\in\N$ and $w\in\Z^+$ such that $x^4+16y^3z+64yz^3$ is a fourth power.
\endproclaim
\Remark\ 1.4. Theorem 1.4(i) implies that $xy,\ 2xy$ and $(x^2+y^2)(x^2+z^2)$ are all suitable.
It seems that (1.4) does hold. We conjecture that $P(x,y,z)$ in Theorem 1.4(i) may be replaced by any of the following polynomials
$$\aligned& (x-y)(x+y-3z),\ (x-y)(x+2y-z),\ (x-y)(x+2y-2z),
\\ &(x-y)(x+2y-7z),\ (x-y)(x+3y-3z),\ (x-y)(x+4y-6z),
\\ &(x-y)(x+5y-2z),\ (x-2y)(x+2y-z),\ (x-2y)(x+2y-2z),
\\ &(x-2y)(x+3y-3z),\ (x+y-z)(x+2y-2z),\ (x-y)(x+y+3z-3w),
\\&(x-y)(x+3y-z-5w),\ (x-y)(3x+3y-3z-5w),
\\&(x-y)(3x+5y-3z-7w),\ (x-y)(3x+7y-3z-9w).
\endaligned\tag1.5$$
In contrast with Theorem 1.4(v), we also conjecture that $x^2y^2+9y^2z^2+9z^2x^2$ is suitable (cf. [S16, A268507]). See [S16, A273110 and A272351] for some data related to
parts (iv) and (vi) of Theorem 1.4.

\proclaim{Theorem 1.5} {\rm (i)} Any $n\in\N$ can be written as $x^2+y^2+z^2+w^2$ with $x,y,z,w\in\N$ such that $x^2-y^2$ is an even square.
Also, $2(x^2-y^2)$, $3(x^2-y^2)$, $x^2-3y^2$ and $3x^2-2y^2$ are all suitable.

{\rm (ii)} All the polynomials
$$\gather x^2+2y^2,\ x^2+3y^2,\ x^2+5y^2,\ x^2+6y^2,\ x^2+8y^2,\ x^2+12y^2,\ 2x^2+7y^2,
\\  3x^2+4y^2,\ 4x^2+5y^2,\ 4x^2+9y^2,\ 5x^2+11y^2,\ 6x^2+10y^2,\ 7x^2+9y^2
\endgather$$
are suitable.
\endproclaim

We will prove Theorems 1.1-1.3 and Theorems 1.4-1.5 in Sections 2 and 3 respectively.
Section 4 contains some open conjectures for further research.

\heading{2. Proofs of Theorems 1.1-1.3}\endheading

\medskip
\noindent {\it Proof of Theorem} 1.1. For $n=0,1,2,\ldots,4^{m/\gcd(2,m)}-1$, the desired result can be verified directly.

Now let $n\gs 4^{m/\gcd(2,m)}$ be an integer and assume that the desired result holds for smaller values of $n$.
\medskip

{\it Case}\ 1. $4^{m/\gcd(2,m)}\mid n$.

By the induction hypothesis, we can write
$$\f n{4^{m/\gcd(2,m)}}=ax^{m}+y^2+z^2+w^2\ \ \t{with}\ x,y,z,w\in\N.$$
It follows that
$$n = a\l(2^{2/\gcd(2,m)}x\r)^m+\l(2^{m/\gcd(2,m)}y\r)^2+\l(2^{m/\gcd(2,m)}z\r)^2+\l(2^{m/\gcd(2,m)}w\r)^2.$$

{\it Case}\ 2. $4^{m/\gcd(2,m)}\nmid n$ and $n\not=4^k(8l+7)$ for any $k,l\in\N$.

In this case, $n\not\in E(1,1,1)$ and hence there are $y,z,w\in\N$ such that $n=a\times0^m+y^2+z^2+w^2$.
\medskip

{\it Case} 3. $n=4^k(8l+7)$ with $k,l\in\N$ and $k<m/\gcd(2,m)$.

If $k\in\{0,1\}$, or $k=2$ and $m>4=a$, then $n-a\not\in E(1,1,1)$ by (1.2), and hence $n=a\times1^{m}+y^2+z^2+w^2$ for some $y,z,w\in\N$.
When $a=1$, $k=2$ and $m\in\{5,6\}$, we have $n-2^m= 4^2(8l+7-2^{m-4})\not\in E(1,1,1)$, and hence $n=a2^m+y^2+z^2+w^2$ for some $y,z,w\in\N$.
If $k\in\{3,4\}$ and $m=5$, then  $n-a2^m\not\in E(1,1,1)$ by (1.2), and hence $n=a2^m+y^2+z^2+w^2$ for some $y,z,w\in\N$.

In view of the above, we have completed our induction proof of Theorem 1.1. \qed

\proclaim{Lemma 2.1} We have
$$E(1,1,2)=\{4^k(16l+14):\ k,l\in\N\}.\tag2.1$$
\endproclaim
\Remark\ 2.1. (2.1) can be found in L. E. Dickson [D39, pp.\,112-113].

\medskip
\noindent {\it Proof of Theorem} 1.2. (i) Let $a\in\{1,2\}$.
We claim that any $n\in\N$ can be written as $x^2+y^2+z^2+w^2$ with $x,y,z,w\in\N$ such that $a(x-y)$ is a square, and want to prove this by induction.

For every $n=0,1,\ldots,15$, we can verify the claim directly.

Now we fix an integer $n\gs16$ and assume that the claim holds for smaller values of $n$.
\medskip

{\it Case} 1. $16\mid n$.

In this case, by the induction hypothesis, there are $x,y,z,w\in\N$ with $a(x-y)$ a square such that
$n/16=x^2+y^2+z^2+w^2$, and hence $n=(4x)^2+(4y)^2+(4z)^2+(4w)^2$ with $a(4x-4y)$ a square.
\medskip

{\it Case} 2. $16\nmid n$ and $n\not\in E(1,1,2)$.

In this case, there are $x,y,z,w\in\N$ with $x=y$ and $n=x^2+y^2+z^2+w^2$, thus $a(x-y)=0^2$ is a square.
\medskip

{\it Case} 3. $16\nmid n$ and $n\in E(1,1,2)$.

In this case, $n=4^k(16l+14)$ for some $k\in\{0,1\}$ and $l\in\N$. Note that $n/2-(2/a)^2\not\in E(1,1,1)$ by (1.2).
So, $n/2-(2/a)^2=t^2+u^2+v^2$ for some $t,u,v\in\N$ with $t\gs u\gs v$. As $n/2-(2/a)^2\gs 8-4>3$, we have $t\gs2\gs 2/a$.
Thus
$$n=2\l(\l(\f2a\r)^2+t^2\r)+2(u^2+v^2)=\l(t+\f2a\r)^2+\l(t-\f2a\r)^2+(u+v)^2+(u-v)^2$$
with $a((t+2/a)-(t-2/a))=2^2$.

So far we have proved part (i) of Theorem 1.2.

(ii) We can easily verify the desired result in Theorem 1.2(ii) for all $n=0,1,\ldots,63$.

Now let $n\gs 64$ and assume that any $r=0,1,\ldots,n-1$ can be written as $x^2+y^2+z^2+w^2\ (x,y,z,w\in\Z)$ with $x+y\in\{ct^3:\ t\in\Z\}$.
If $64\mid n$, then $n/64$ can be written as $x^2+y^2+z^2+w^2\ (x,y,z,w\in\Z)$ with $x+y=ct^3$ for some $t\in\Z$, hence
$n=(8x)^2+(8y)^2+(8z)^2+(8w)^2$ with $8x+8y=c(2t)^3$.

Now we consider the case $64\nmid n$. We claim that $\{2n,2n-c^2,2n-64c^2\}\not\se E(1,1,1)$.
If $\{2n,2n-1\}\se E(1,1,1)$, then (1.2) implies that $2n=4^k(8l+7)$ for some $k\in\{2,3\}$ and $l\in\N$, hence
$2n-64$ is $4^2(8l+3)$ or $4^3(8l+6)$, and thus $2n-64\not\in E(1,1,1)$.
If $2n=4^k(8l+7)$ for some $k\in\{1,2\}$ and $l\in\N$, then $2n-4$ is $4(8l+6)$ or $4(8(4l+3)+3)$, and hence $2n-2^2\not\in E(1,1,1)$;
if $2n=4^3(8l+7)$ with $l\in\N$, then $2n-64\times2^2=4^3(8l+3)\not\in E(1,1,1)$.
When $2n=4^k(8l+7)$ for some $k\in\{1,2,3\}$ and $l\in\N$, we have $2n-4^2\not\in E(1,1,1)$ since
$$2n-16= 4^k(8l+7)-16=\cases 4(8l+7)-16=4(8l+3)&\t{if}\ k=1,
\\4^2(8l+7)-16=4^2(8l+6)&\t{if}\ k=2,
\\4^3(8l+7)-16=4^2(8(4l+3)+3)&\t{if}\ k=3.
\endcases$$
So the claim is true.

By the claim, for some $\da\in\{0,1,8\}$, we can write
$2n-\da^2c^2$ as the sum of three squares two of which have the same parity. Hence we may write $2n-\da^2c^2=(2x-\da c)^2+y^2+z^2$ with $x,y,z\in\Z$ and $y\eq z\pmod2$.
It follows that
$$\align n=&\f{(2x-\da c)^2+\da^2c^2}2+\f{y^2+z^2}2
\\=&x^2+(\da c-x)^2+\l(\f{y+z}2\r)^2+\l(\f{y-z}2\r)^2
\endalign$$
with $x+(\da c-x)=c\da\in\{ct^3:\ t=0,1,2\}$. This concludes the induction step.

(iii) For $n=0,1,\ldots,4^m-1$, we can verify the desired result in Theorem 1.2(iii) directly via a computer.

Fix $n\in\N$ with $n\gs 4^m$, and assume the required result for smaller values of $n$.

If $4^m\mid n$, then by the induction hypothesis $n/4^m$ can be written as $x^2+y^2+z^2+w^2$ with $x,y,z,w\in\Z$ such that $x+2y=dt^m$ for some $t\in\Z$,
and hence $n=(2^mx)^2+(2^my)^2+(2^mz)^2+(2^mw)^2$ with $2^mx+2(2^my)=2^m(x+2y)=d(2t)^m$.

Now we assume $4^m\nmid n$. In light of (1.2), $\{5n,5n-d^2,5n-4^md^2\}\not\se E(1,1,1)$. In fact, for any $l\in\N$ neither $8l+7-d^2$
nor $4(8l+7)-d^2$ belongs to $E(1,1,1)$; if $5n=4^2(8l+7)$ with $l\in\N$ then $5n-2^2=4(8(4l+3)+3)\not\in E(1,1,1)$ and $5n-4^31^2=4^2(8l+3)\not\in E(1,1,1)$.
So, for some $\da\in\{0,1,2^m\}$ we can write $5n-\da^2d^2$ as $x^2+y^2+z^2$ with $x,y,z\in\Z$.
Note that a square is congruent to one of $0,1,-1$ modulo $5$.
It is easy to see that one of $x^2,y^2,z^2$ is congruent to $-(\da d)^2$ modulo $5$.
Without loss of generality, we may assume that $(\da d)^2+x^2\eq y^2+z^2\eq0\pmod 5$.
As $x^2\eq(2\da d)^2\pmod 5$ and $y^2\eq(2z)^2\pmod 5$, without loss of generality we simply assume that $x\eq 2\da d\pmod 5$
(otherwise we use $-x$ instead of $x$) and $y\eq 2z\pmod 5$.
Thus $r=(2x+\da d)/5,\ s=(2\da d-x)/5$, $u=(2y+z)/5$ and $v=(2z-y)/5$ are all integers.
Note that
$$r^2+s^2+u^2+v^2=\f{(\da d)^2+x^2}5+\f{y^2+z^2}5=n \ \ \t{with}\ \  r+2s=\da d\in\{dt^m:\ t=0,1,2\}.$$
This concludes our proof of the third part of Theorem 1.2.

(iv) For each $n=0,1,2$, there are $x,z\in\{0,1\}$ such that $n=x^2+0^2+z^2+0^2$ with $(x+z)(0+0)=0^2$. Note also that
$$3=1^2+1^2+(-1)^2+0^2\quad\t{with}\ (1+(-1))(1+0)=0^2.$$

Now let $n\in\{4,5,6,\ldots\}$ and assume that any $r=0,\ldots,n-1$ can be written as $x^2+y^2+z^2+w^2$ with $x,y,z,w\in\Z$ such that
$(x+z)(y+w)$ is a square.

If $4\mid n$, then $n/4$ can be written as $x^2+y^2+z^2+w^2\ (x,y,z,w\in\Z)$ with $(x+z)(y+w)$ a square,
and hence $n=(2x)^2+(2y)^2+(2z)^2+(2w)^2$ with $(2x+2z)(2y+2w)=2^2(x+z)(y+w)$ a square.

If $n$ is odd, then by (2.1) there are $x,y,z,w\in\Z$ with $w=-y$ such that $n=x^2+2y^2+z^2=x^2+y^2+z^2+w^2$ with $(x+z)(y+w)=0^2$.

Now we consider the case $n=2m$ with $m$ odd. By (2.1), we can write $m$ as $x^2+y^2+z^2+w^2$ with $x,y,z,w\in\Z$ and $w=y$. Therefore,
$$n=2m=(x+y)^2+(x-y)^2+(z+w)^2+(z-w)^2=(x+y)^2+(z+w)^2+(y-x)^2+(w-z)^2$$
with
$$((x+y)+(y-x))((z+w)+(w-z))=2y(2w)=(2y)^2.$$
This ends our induction proof of Theorem 1.2(iv).

(v) Clearly,
$$1=4^0(1+4\times0^2+0^2)+0^2,\ 2=4^0(1+4\times0^2+1^2)+0^2,\ 3=4^0(1+4\times0^2+1^2)+1^2.$$

Now let $n\in\Z^+$ with $n\gs 4$. If $4\mid n$ and $n/4=4^k(1+4x^2+y^2)+z^2$ for some $k,x,y,z\in\N$, then $n=4^{k+1}(1+4x^2+y^2)+(2z)^2$.
If $n\eq2,3\pmod 4$, then $n-1\eq1,2\pmod 4$ and hence for some $x,y,z\in\Z$ we have $n-1=(2x)^2+y^2+z^2$ and thus $n=4^0(1+4x^2+y^2)+z^2$.
By Dickson [D39, pp.\,112--113], if $q\in\N$ is congruent to $1$ modulo $4$ then $q\not\in E(1,4,16)$.
Thus, when $n\eq1\pmod4$, we have $n-4\not\in E(1,4,16)$, hence there are $x,y,z\in\Z$ such that $n-4=16x^2+4y^2+z^2$, i.e., $n=4(1+4x^2+y^2)+z^2$.
This proves Theorem 1.2(v) by induction.

\medskip

In view of the above, we have completed the proof of Theorem 1.2. \qed
\medskip

\Remark\ 2.2. In contrast with (2.1), we conjecture that any odd integer $n>1248$ can be written as
$x^2+y^2+2z^2$ with $x,y,z\in\N$ and $x\gs y\gs z$. Modifying the proof of Theorem 1.2(iii) slightly, we see that
under this conjecture the polynomial $(x-z)(y-w)$ is suitable and also any $n\in\N$
can be written as $x^2+y^2+z^2+w^2$ with $w\in\Z$, $x,y,z\in\N$ and $x\ls y\gs z\gs|w|$ such that $(x+y)(z+w)$ is a square.

\proclaim{Lemma 2.2} We have
$$\align E(1,2,6)=&\{4^k(8l+5):\ k,l\in\N\},\tag2.2
\\ E(2,3,6)=&\{3q+1:\ q\in\N\}\cup\{4^k(8l+7):\ k,l\in\N\},\tag2.3
\\E(1,2,3)=&\{4^k(16l+10):\ k,l\in\N\},\tag2.4
\\ E(1,3,6)=&\{3q+2:\ q\in\N\}\cup\{4^k(16l+14):\ k,l\in\N\},\tag2.5
\\E(1,5,5)=&\{n\in\N:\ n\eq2,3\ (\mo\ 5)\}\cup\{4^k(8l+7):\ k,l\in\N\}.\tag2.6
\endalign$$
\endproclaim
\Remark\ 2.3. (2.2)-(2.6) are known results, see, e.g., Dickson [D39, pp.\,112-113].

\medskip
\noindent {\it Proof of Theorem} 1.3. (i) For $n=0,\ldots,4^m-1$ we can easily verify the desired result in Theorem 1.3(i) directly.

Now let $n\in\N$ with $n\gs 4^m$. Assume that any $r\in\{0,\ldots,n-1\}$ can be written as $x^2+y^2+z^2+w^2$
with $x,y,z,w\in\Z$ such that $x+y+cz\in\{dt^m:\ t\in\Z\}$. If $4^m\mid n$, then there are $x,y,z,w\in\Z$
with $x^2+y^2+z^2+w^2=n/4^m$ such that $x+y+cz=dt^m$ for some $t\in\Z$, and hence
$$n=(2^mx)^2+(2^my)^2+(2^mz)^2+(2^mw)^2$$
with $2^mx+2^my+c(2^mz)=2^m(x+y+cz)=d(2t)^m$.
Below we suppose that $4^m\nmid n$.
\medskip

{\it Case}\ 1. $c=1$.

In this case, it suffices to show that there are $x,y,z\in\Z$ and $\da\in\{0,1,2^m\}$ such that
$$n=x^2+(y+z)^2+(z-y)^2+(\da d-2z)^2=x^2+2y^2+6z^2-4\da dz+\da^2d^2.\tag2.7$$
(Note that $(y+z)+(z-y)+(\da d-2z)=\da d\in\{dt^m:\ t\in\Z\}$.) Suppose that this fails for $\da=0$. By (2.2),
$n=4^k(8l+5)$ for some $k,l\in\N$ with $k<m$.

We first handle the subcase $d=1$. Clearly,
$$3n-1=\cases3(8l+5)-1=2(12l+7)&\t{if}\ k=0,\\3\times4(8l+5)-1 = 8(12l+7)+3&\t{if}\ k=1.\endcases$$
Thus, if $k\in\{0,1\}$, then $3n-1\not\in E(2,3,6)$ by (2.3),
hence for some $x,y,z\in\Z$ we have
$$3n-1=3x^2+6y^2+2(3z-1)^2=3(x^2+2y^2+2(3z^2-2z))+2$$
and thus
$$n=x^2+2y^2+6z^2-4z+1=x^2+(y+z)^2+(z-y)^2+(1-2z)^2$$
which gives (2.7) with $\da=1$.
When $k=2$ and $m=3$, we have
$$3n-64=3\times4^2(8l+5)-64=4^2(8(3l+1)+3)\not\in E(2,3,6)$$
in view of (2.3), hence for some $x,y,z\in\Z$ we have
$$3n-4^3=3x^2+6y^2+2(3z-8)^2=3(x^2+2y^2+2(3z^2-16z))+2\times 4^3$$
and thus
$$n=x^2+2y^2+6z^2-32z+64=x^2+(y+z)^2+(z-y)^2+(2^3-2z)^2$$
which gives (2.7) with $\da=2^m$.

Now we handle the subcase $d=2$. Clearly,
$$3n-4=3\times 4^k(8l+5)-4=\cases 8(3l+1)+3&\t{if}\ k=0,
\\4(8(3l+1)+6)&\t{if}\ k=1,
\\4(8(12l+7)+3)&\t{if}\ k=2.
\endcases$$
It follows that $3n-4\not\in E(2,3,6)$ by (2.3). So, for some $x,y,z\in\Z$ we have
$$3n-4=3x^2+6y^2+2(3z-2)^2=3x^2+6y^2+18z^2-24z+8$$
and hence
$$n=x^2+2y^2+6z^2-8z+4=x^2+(y+z)^2+(z-y)^2+(2-2z)^2$$
with $(y+z)+(z-y)+(2-2z)=2\times1^m$.
\medskip

{\it Case}\ 2. $c=2$.

If $n\not\in E(1,2,3)$, then for some $x,y,z\in\Z$ we have
$$n=x^2+2y^2+3z^2=x^2+(y+z)^2+(z-y)^2+(-z)^2$$
with $(y+z)+(z-y)+2(-z)=d\times0^m$. Now let $n\in E(1,2,3)$. By (2.4), $n=4^k(16l+10)$ for some $k,l\in\N$ with $k<m$.
\medskip

{\it Subcase} 2.1. $m=2$ and $d\in\{3,6\}$.

For $n=16,\ldots,23$ we can verify the desired result directly. Now let $n\gs 24=216/9$.
No matter $k=0$ or $k=1$, we have
$$n-\f{216}{d^2}=4^k(16l+10)-\f{216}{d^2}\not\in E(1,2,3)$$
by (2.4). So, there are $x,y,z\in\Z$ such that
$$n-\f{216}{d^2}=x^2+2y^2+3\l(\f{6}d-z\r)^2$$
and hence
$$n=x^2+2y^2+3z^2-\f{36}dz+\l(\f{18}d\r)^2=x^2+(y+z)^2+(z-y)^2+\l(\f{18}d-z\r)^2$$
with $(y+z)+(z-y)+2(18/d-z)=d(6/d)^2$.

\medskip

{\it Subcase} 2.2. $d=1$.

If $k=0$, then $6n-1=6(16l+10)-1\eq3\pmod 8$ and hence $6n-1\not\in E(2,3,6)$ by (2.3).
When $k=1$, we have
$$6n-4^m=6\times4(16l+10)-4^m=4^2(8(3l+1)+7-4^{m-2})\not\in E(2,3,6)$$
by (2.3). If $k=2$ and $m=3$, then
$$6n-4^m=6\times 4^2(16l+10)-4^3=4^3(8(3l+1)+6)\not\in E(2,3,6)$$
by (2.3). So, for some $\da\in\{1,2^m\}$ we have
$6n-\da^2\not\in E(2,3,6)$. Hence there are $x,u,v\in\Z$ such that
$6n-\da^2=6x^2+3u^2+2v^2$. As $u\eq\da\pmod 2$, we can write $u=2y+\da$ with $y\in\Z$.
Since $v^2\eq\da^2\pmod3$, without loss of generality we may assume that $v=3z+\da$ with $z\in\Z$.
Therefore,
$$6n=6x^2+3(2y+\da)^2+2(3z+\da)^2+\da^2$$
and hence
$$n=x^2+2y^2+3z^2+2\da y+2\da z+\da^2=x^2+(y+z+\da)^2+(z-y)^2+(-z)^2$$
with $(y+z+\da)+(z-y)+2(-z)=\da\in\{t^m:\ t=1,2\}$.
\medskip

{\it Subcase} 2.3. $d=2$.

Observe that
$$3n-2=3\times 4^k(16l+10)-2=\cases 4(12l+7)&\t{if}\ k=0,
\\16(12l+7)+6&\t{if}\ k=1.
\endcases$$
If $k=2$ and $m=3$, then
$$3n-2^7=3\times4^2(16l+10)-2^7=4^2(16(3l+1)+6).$$
Combining this with (2.5), we find that for some $\da\in\{1,2^m\}$ we have $3n-2\da^2\not\in E(1,3,6)$.
Thus, there are $x,y,z\in\Z$ for which
$$3n-2\da^2=3x^2+6y^2+(3z-\da)^2=3x^2+6y^2+9z^2-6\da z+\da^2$$
and hence
$$n=x^2+2y^2+3z^2-2\da z+\da^2=x^2+(y+z)^2+(z-y)^2+(\da-z)^2$$
with $(y+z)+(z-y)+2(\da-z)=2\da\in\{2t^m:\ t=1,2\}$.
\medskip

{\it Subcase} 2.4. $d=4$.

Observe that
$$3n-8=3\times 4^k(16l+10)-8=\cases 16(3l+1)+6&\t{if}\ k=0,
\\4^2(12l+7)&\t{if}\ k=1,\\4(16(12l+7)+6)&\t{if}\ k=2.
\endcases$$
Clearly, $3n-8\not\in E(1,3,6)$ by (2.5).
Thus, for some $x,y,z\in\Z$ we have
$$3n-8=3x^2+6y^2+(3z-2)^2=3x^2+6y^2+9z^2-12z+4$$
and hence
$$n=x^2+2y^2+3z^2-4z+4=x^2+(y+z)^2+(z-y)^2+(2-z)^2$$
with $(y+z)+(z-y)+2(2-z)=4\times1^m$.

Combining the above, we have proved part (i) of Theorem 1.3.

(ii) Suppose that any odd integer $n>2719$ can be represented by $x^2+y^2+10z^2$ with $x,y,z\in\Z$.
We want to prove by induction the claim that each $n\in\N$ can be written as $x^2+y^2+z^2+w^2$ ($x,y,z,w\in\Z$) with $x+3y$ a square.

 For $n=0,1,\ldots,15$, the claim can be verified via a computer.

 Now fix an integer $n\gs16$ and assume that the claim holds for smaller values of $n$.

 If $16\mid n$, then by the induction hypothesis there are $x,y,z,w\in\Z$ with $n/16=x^2+y^2+z^2+w^2$ such that $x+3y$ is a square,
 and hence $n=(4x)^2+(4y)^2+(4z)^2+(4w)^2$ with $4x+3(4y)=4(x+3y)$ a square.

 Now we let $16\nmid n$. If $2\nmid n$ and $n\ls 2719$, then we can easily verify that $5n$ or $5n-8$ can be written as $2x^2+5y^2+5z^2$ with $x,y,z\in\Z$.
 If $2\nmid n$ and $n>2719$, then there are $x,y,z\in\Z$ such that $n=10x^2+y^2+z^2$ and hence $5n=2(5x)^2+5y^2+5z^2$.
 If $n$ is even and $n$ is not of the form $4^k(16l+6)$ $(k,l\in\N)$, then by Dickson [D27] there are $x,y,z\in\Z$ such that $n=10x^2+y^2+z^2$ and hence $5n=2(5x)^2+5y^2+5z^2$.
 When $n=4^k(16l+6)$ for some $k\in\{0,1\}$ and $l\in\N$, clearly
 $$\f{5n-8}2=5\times 4^k(8l+3)-4\not\in E(1,5,5)$$
 by (2.6), thus there are $x,y,z\in\Z$
 such that $(5n-8)/2=x^2+5y^2+5z^2$ and hence $5n-8=2x^2+5(y+z)^2+5(y-z)^2$.

 Since $5n$ or $5n-8$ can be written as $2x^2+5y^2+5z^2$ with $x,y,z\in\Z$, for some $\da\in\{0,2\}$ and $x,y,z\in\Z$ we have
 $$10n-\da^4=2(2x^2+5y^2+5z^2)=(2x)^2+10y^2+10z^2.$$
 As $(2x)^2\eq -\da^4\eq(3\da^2)^2\pmod{10}$, without loss of generality we may assume that $2x=10w+3\da^2$ with $w\in\Z$.
 Then
 $$10n=\da^4+(10w+3\da^2)^2+10y^2+10z^2$$
 and hence
 $$n=10w^2+y^2+z^2+6\da^2w+\da^4=(3w+\da^2)^2+(-w)^2+y^2+z^2$$
 with $(3w+\da^2)+3(-w)=\da^2$ a square.
This concludes the induction step.

(iii) Suppose that any integer $n\gs1190$ not divisible by $16$ can be written as $x^2+10y^2+(2z^2+5^3r^4)/7$ with $r\in\{0,1,2,3\}$ and $x,y,z\in\Z$.
 We want to prove by induction the claim that each $n\in\N$ can be written as $x^2+y^2+z^2+w^2$ with $x,y,z,w\in\Z$ such that $x+3y+5z$ is a square.

 For $n=0,1,\ldots,1189$, the claim can be verified via a computer.

 Now fix an integer $n\gs1190$ and assume that the claim holds for smaller values of $n$.

 If $16\mid n$, then by the induction hypothesis there are $x,y,z,w\in\Z$ with $n/16=x^2+y^2+z^2+w^2$ such that $x+3y+5z$ is a square,
 and hence $n=(4x)^2+(4y)^2+(4z)^2+(4w)^2$ with $4x+3(4y)+5(4z)=4(x+3y+5z)$ a square.

 Below we let $16\nmid n$. Then, for some $r\in\{0,1,2,3\}$ and $x,y,z\in\Z$ we have
 $$n=x^2+10y^2+\f{2z^2+5^3r^4}7,\ \t{i.e.},\ 7n-5(5r^2)^2=7x^2+70y^2+2z^2.$$
 Thus,
 $$7(n-(5r^2)^2)=7x^2+70y^2+2\l(z^2-(5r^2)^2\r).$$
 As $5r^2$ is congruent to $z$ or $-z$ modulo $7$, without loss of generality we may assume that $z=7t-5r^2$ for some $t\in\Z$.
 It follows that
 $$n-(5r^2)^2=x^2+10y^2+\f27\l((7t-5r^2)^2-(5r^2)^2\r)=x^2+10y^2+14t^2-4(5r^2)t$$
and hence
$$n=x^2+10y^2+10t^2+(2t-5r^2)^2=x^2+(3y+t)^2+(3t-y)^2+(5r^2-2t)^2$$
with
$$(3y+t)+3(3t-y)+5(5r^2-2t)=(5r)^2.$$
This concludes the induction step.

The proof of Theorem 1.3 is now complete. \qed
\medskip

\Remark\ 2.4. Those $3z^2-2z$ with $z\in\Z$ appeared in the proof of Theorem 1.3(i) are called generalized octagonal numbers, one may consult [S15] and [S16a] for related results.

\heading{3. Proofs of Theorems 1.4-1.5}\endheading

\proclaim{Lemma 3.1} We have
$$ E(1,1,5)=\{4^k(8l+3):\ k,l\in\N\},\tag3.1$$
and
$$E(1,1,10)\cap2\Z=\{4^k(16l+6):\ k,l\in\N\}.\tag3.2$$
\endproclaim
\Remark\ 3.1. (3.1) can be found in Dickson [D39, pp.\,112-113], and (3.2) was a conjecture of Ramanujan [R] proved by Dickson [D27].
\medskip

Combining (1.2), (2.1), (2.2), (2.4) and (3.1)-(3.2), we immediately get the following lemma.

\proclaim{Lemma 3.2} The six sets
$$E(1,1,1),\ E(1,1,2),\ E(1,2,3),\ E(1,2,6),\ E(1,1,5)\ \ \t{and}\ \ E(1,1,10)\cap 2\Z\tag3.3$$
are pairwise disjoint.
\endproclaim

\proclaim{Lemma 3.3} Let $n\in\N$. Then $n\not\in E(1,2,6)$ if and only if $n=x^2+y^2+z^2+w^2$ for some $x,y,z,w\in\N$ with $x+y=z$.
Also, $n\not\in E(1,2,3)$ if and only if $n=x^2+y^2+z^2+w^2$ for some $x,y,z,w\in\Z$ with $x+y=2z$.
\endproclaim
\Proof. (i) Assume that $n\not\in E(1,2,6)$. Then, there are $x,y,z\in\N$ for which
$$n=x^2+2y^2+6z^2=x^2+(y+z)^2+|y-z|^2+(2z)^2.$$
Clearly $(y+z)+|y-z|=2z$ if $y\ls z$, and $|y-z|+2z=y+z$ if $y>z$. Therefore $n=x^2+u^2+v^2+w^2$
for some $u,v,w\in\N$ with $u+v=w$.

Now suppose that $n=x^2+y^2+z^2+w^2$ with $x,y,z,w\in\N$ and $x+y=z$. If $x\eq y\pmod 2$, then
$$n=2\l(\f{x+y}2\r)^2+2\l(\f{x-y}2\r)^2+(x+y)^2+w^2=w^2+2\l(\f{x-y}2\r)^2+6\l(\f{x+y}2\r)^2$$
and hence $n\not\in E(1,2,6)$. When $x\not\eq y\pmod 2$, without loss of generality we may assume that $y\eq z\pmod 2$, hence
$$n=(z-y)^2+2\l(\f{y+z}2\r)^2+2\l(\f{y-z}2\r)^2+w^2=w^2+2\l(\f{y+z}2\r)^2+6\l(\f{y-z}2\r)^2$$
and thus $n\not\in E(1,2,6)$.

(ii) If $n\not\in E(1,2,3)$, then there are $x,y,z\in\Z$ for which
$$n=x^2+2y^2+3z^2=x^2+(y+z)^2+(z-y)^2+z^2$$
with $(y+z)+(z-y)=2z$.

If $n=x^2+y^2+z^2+w^2$ with $x,y,z,w\in\Z$ and $x+y=2z$, then
$$n=2\l(\f{x+y}2\r)^2+2\l(\f{x-y}2\r)^2+z^2+w^2=w^2+2\l(\f{x-y}2\r)^2+3z^2$$
and hence $n\not\in E(1,2,3)$.

In view of the above, we have completed the proof of Lemma 3.3. \qed

\medskip
\noindent {\it Proof of Theorem} 1.4. (i) Clearly, $5y^2=(2y)^2+y^2$ and $10y^2=(3y)^2+y^2$.
In view of Lemmas 3.2 and 3.3, each $n\in\N$ can be written as
$x^2+y^2+z^2+w^2$ with $x,y,z,w\in\N$ such that $P(x,y,z)=0$, where $P(x,y,z)$ may be any of the polynomials listed in (1.3).
For example, for $P(x,y,z)=(x-2y)(x+y-z)$ we use the fact that $E(1,1,5)\cap E(1,2,6)=\em$.

(ii) Suppose that (1.4) holds.
If $n=4^k(8l+5)$ for some $k,l\in\N$, then there are $x,y,z\in\N$ such that $8l+5=x^2+y^2+13z^2$
and hence
$$n = (2^kx)^2+(2^ky)^2+13(2^kz)^2 = (2^kx)^2+(2^ky)^2+(2^k\times 3z)^2+(2^{k+1}z)^2$$
with $2(2^k\times 3z)=3(2^{k+1}z)$. If $n\in\N\sm\{4^k(8l+5):\ k,l\in\N\}$, then $n\not\in E(1,2,6)$ by (2.2), and hence by Lemma 3.3 there are $x,y,z,w\in\N$ such that
$n=x^2+y^2+z^2+w^2$ with $x+y=z$. Therefore, any $n\in\N$ can be written as $x^2+y^2+z^2+w^2$ ($x,y,z,w\in\N)$ with $(2x-3y)(x+y-z)=0$.

(iii) If $n=2x^2$ for some $x\in\Z^+$, then $n=x^2+0^2+x^2+0^2$ with $x\in\Z^+$ and $(x-0)x=x^2$.
If $n\not\in E(1,1,2)$ and $n\not=2x^2$ for all $x\in\N$, then there are $x,z,w\in\N$ with $z>0$ such that
$n=x^2+x^2+z^2+w^2$ with $(x-x)z=0^2$.

Now suppose that $n\in E(1,1,2)$. Then $n\not\in E(1,2,6)$ by Lemma 3.2, and hence by Lemma 3.3 there are $x,y,z,w\in\N$ with $x=y+z$
such that $x^2+y^2+z^2+w^2=n$. If $y$ and $z$ are both zero, then $n=w^2\not\in E(1,1,2)$ which leads a contradiction.
Without loss of generality we assume that $z>0$. Note that $(x-y)z=z^2$ is a square. This concludes the proof of Theorem 1.4(iii).

(iv) If $n\not\in E(1,1,1)$, then there are $y,z,w\in\N$ with $y>0$ such that $n=0^2+y^2+z^2+w^2$. As
$$(0+4y+4z)^2+(9\times0+3y+3z)^2=(5y+5z)^2,$$
by Pythagoras' theorem there is a right triangle with positive integer sides $4y+4z$, $3y+3z$ and $5y+5z$.

Now let $n\in E(1,1,1)$. Then $n\not\in E(1,2,6)$ by Lemma 3.2, and hence by Lemma 3.3 there are $x,y,z,w\in\N$ with $x=y+z$ such that
$n=x^2+y^2+z^2+w^2$. Clearly $y>0$ since $n\in E(1,1,1)$.
Observe that
$$(x+4y+4z)^2+(9x+3y+3z)^2=(5x)^2+(12x)^2=(13x)^2.$$
So $x+4y+4z$ and $9x+3y+3z$ are the two legs of a right triangle with positive integer sides.

(v) Fix $n\in\Z^+$. If $n\not\in E(1,1,1)$, then there are $x,y\in\N$ and $w\in\Z^+$ such that $n=x^2+y^2+0^2+w^2$
and hence $x^2y^2+y^20^2+0^2x^2=x^2y^2+4y^20^2+4\times0^2x^2=(xy)^2$ is a square.

Now suppose that $n\in E(1,1,1)$. Then $n\not\in E(1,2,6)$ by Lemma 3.2, and hence by Lemma 3.3 there are $x,y,z,w\in\N$ with $x+y=z$ such that
$n=x^2+y^2+z^2+w^2$. Clearly $w>0$ since $n\in E(1,1,1)$.
Observe that
$$\align x^2y^2+y^2z^2+z^2x^2=&(xy)^2+(x^2+y^2)(x+y)^2
\\=& (xy)^2+(x^2+xy+y^2-xy)(x^2+xy+y^2+xy)
\\=& (x^2+xy+y^2)^2.
\endalign$$

As $E(1,2,3)\cap E(1,1,1)=\em$ by Lemma 3.2, we have $n\not\in E(1,2,3)$, and hence by Lemma 3.3 we can write $n$ as $x^2+y^2+z^2+w^2\ (x,y,z,w\in\Z)$ with $x+y=2z$.
Clearly $w\not=0$ since $n\in E(1,1,1)$. Note that
$$\align x^2y^2+4y^2z^2+4z^2x^2=&(xy)^2+(x^2+y^2)(2z)^2
\\=&(xy)^2+(x^2+y^2)(x+y)^2=(x^2+xy+y^2)^2.
\endalign$$
This concludes the proof of Theorem 1.4(v).

(vi) If $n\not\in E(1,1,1)$, then there are $x,z,w\in\N$ with $w>0$ such that $n=x^2+0^2+z^2+w^2$ and hence
$x^4+8\times0^3z+8\times0z^3=x^4+16\times0^3z+16\times0z^3 = x^4$ is a fourth power.
Below we assume that $n\in E(1,1,1)$.

As $E(1,1,1)\cap E(1,2,6)=\em$ (by Lemma 3.2), we have $n\not\in E(1,2,6)$. In view of Lemma 3.3, there are $w,x,y,z\in\N$ with $x+z=y$ such that
$n=x^2+y^2+z^2+w^2$. As $n\in E(1,1,1)$, we have $wx\not=0$.
Observe that
$$\align x^4+8yz(y^2+z^2)=& x^4+((y+z)^2-(y-z)^2)((y+z)^2+(y-z)^2)
\\=&x^4+(y+z)^4-(y-z)^4=(y+z)^4.
\endalign$$

As $E(1,1,1)\cap E(1,2,3)=\em$ by Lemma 3.2, we have $n\not\in E(1,2,3)$ and hence
there are $v,w,z\in\N$ such that $n=w^2+2v^2+3z^2=w^2+(v+z)^2+(v-z)^2+z^2$.
Clearly $w>0$ since $n\not\in E(1,1,1)$. Let $y=v+z$ and $x=|v-z|=|y-2z|$. Then
$n=w^2+x^2+y^2+z^2$ and
$$\align x^4+16y^3z+64yz^3=&x^4+4y(2z)\l(2y^2+2(2z)^2\r)
\\=&x^4+\l((y+2z)^2-(y-2z)^2\r)\l((y+2z)^2+(y-2z)^2\r)
\\=&x^4+(y+2z)^4-(y-2z)^4=(y+2z)^4.
\endalign$$

The proof of Theorem 1.4 is now complete. \qed

\medskip
\noindent {\it Proof of Theorem} 1.5. (i) If $n\not\in E(1,1,2)$, then for some $x,y,z,w\in\N$ with $x=y$ we have
$$n =x^2+y^2+z^2+w^2$$
with $a(x^2-y^2)=0^2$ for any $a=1,2,3$, and $3x^2-2y^2=x^2$. If $n\not\in E(1,1,1)$, then $n=x^2+0^2+y^2+z^2$ for some $x,y,z\in\Z$ with $x^2-3\times0^2=x^2$.
If $n\not\in E(1,1,5)$, then there are $x,y,z\in\N$ satisfying $n=5x^2+y^2+z^2=(2x)^2+x^2+y^2+z^2$
with $3((2x)^2-x^2)=(3x)^2$ and $(2x)^2-3x^2=x^2$.
Since
$$E(1,1,2)\cap E(1,1,5)=E(1,1,1)\cap E(1,1,5)=\em$$
by Lemma 3.2, we conclude that both $3(x^2-y^2)$ and $x^2-3y^2$ are suitable polynomials.

Now suppose that $n\in E(1,1,2)$. By (2.1), we have $n=4^k(16l+14)$ for some $k,l\in\N$.
As $n\not\in E(1,1,1)$ and $2\mid n$, there are $u,v,w\in\N$ satisfying $n=(2w)^2+0^2+u^2+v^2$
with $(2w)^2-0^2$ an even square.
Since $n$ is even and $n\not=4^i(16j+6)$ for any $i,j\in\N$, we have $n\not\in E(1,1,10)$ by (3.2).
Thus, there are $x,y,z\in\N$ satisfying $n=10x^2+y^2+z^2=(3x)^2+x^2+y^2+z^2$
with $2((3x)^2-x^2)=(4x)^2$ and $3(3x)^2-2x^2=(5x)^2$.

In view of the above, we have proved part (i) of Theorem 1.5.

(ii) We first prove that
$x^2+ky^2\ (k=2,3,5,6,8,12)$ and $7x^2+9y^2$ are all suitable.

Let $n\in\N$. If $n\not\in E(1,1,1)$, then there are $x,y,z\in\N$
such that $n=x^2+0^2+y^2+z^2$ with $x^2+k0^2=x^2$ for any $k\in\Z^+$, and $9x^2+7\times0^2=(3x)^2$.

Now assume that $n\in E(1,1,1)$. As $n\not\in E(1,1,2)$ by Lemma 3.2, there are $x,y,z,w\in\N$ with $x=y$ and
$n=x^2+y^2+z^2+w^2$, hence $x^2+3y^2=(2x)^2$, $x^2+8y^2=(3x)^2$ and $7x^2+9y^2=(4x)^2$.
As $n\not\in E(1,1,5)$ by Lemma 3.2, there are $x,y,z,w\in\N$ with $2x=y$ and
$n=x^2+y^2+z^2+w^2$, hence
$$x^2+2y^2=(3x)^2,\ y^2+5x^2 = (3x)^2,\ x^2+6y^2= (5x)^2, \ x^2+12y^2=(7x)^2.$$

Next we show that $3x^2+4y^2,4x^2+5y^2$ and $4x^2+9y^2$ are all suitable.

Let $n\in\N$. If $n\not\in E(1,1,1)$, then for some $x,y,z\in\N$ we have $n=x^2+0^2+y^2+z^2$
with $4x^2+b\times0^2=(2x)^2$ for $b=3,5,9$.

Now assume that $n\in E(1,1,1)$. As $n\not\in E(1,1,5)$ by Lemma 3.2, there are $x,y,z,w\in\N$ with $y=2x$ such that
$n=x^2+y^2+z^2+w^2$ and hence $3y^2+4x^2=(4x)^2$ and $4y^2+9x^2=(5x)^2$.
As $n\not\in E(1,1,2)$ by Lemma 3.2, there are $x,y,z,w\in\N$ with $y=x$ such that
$n=x^2+y^2+z^2+w^2$ and hence $4x^2+5y^2=(3x)^2$.

Finally we prove that $2x^2+7y^2$, $5x^2+11y^2$ and $6x^2+10y^2$ are suitable.

Let $n\in\N$. If $n\not\in E(1,1,2)$, then there are $x,y,z,w\in\N$ with $y=x$ such that
$n=x^2+y^2+z^2+w^2$ and hence
$$2x^2+7y^2=(3x)^2\ \ \t{and}\ \  5x^2+11y^2=(4x)^2=6x^2+10y^2.$$

Below we assume that $n\in E(1,1,2)$. Then $2\mid n$ by (2.1). As $n\not\in E(1,1,5)$ by Lemma 3.2, there are $x,y,z,w\in\N$ with $y=2x$ such that
$n=x^2+y^2+z^2+w^2$ and hence $5x^2+11y^2=(7x)^2$.
Since $2\mid n$ and $n\not\in E(1,1,10)\cap 2\Z$ by Lemma 3.2,
there are $x,y,z,w\in\N$ with $y=3x$ such that
$n=x^2+y^2+z^2+w^2$ and hence $2y^2+7x^2=(5x)^2$ and $6y^2+10x^2=(8x)^2$.

In view of the above, we have completed the proof of Theorem 1.5. \qed

\heading{4. Some open conjectures}\endheading

Motivated by our results in Section 1, we pose the following conjectures for further research.

\proclaim{Conjecture 4.1} Let $a,b\in\Z^+$ with  $\gcd(a,b)$ squarefree. Then $ax+by$ is suitable if and only if $\{a,b\}=\{1,2\},\{1,3\},\{1,24\}$.
Also, $ax-by$ is suitable if and only if $(a,b)$ is among the ordered pairs
$$(1,1),\ (2,1), \ (2,2),\ (4,3), \ (6,2).$$
\endproclaim
\Remark\ 4.1. By Theorem 1.2(i), both $x-y$ and $2x-2y$ are suitable. Though we have Theorem 1.2(iii) and Theorem 1.3(ii), we are not able to show that
$x+2y$ or $2x-y$ or $x+3y$ is suitable. See [S16, A273404] for the number of ordered ways to write $n=x^2+y^2+z^2+w^2$ with $x,y,z,w\in\N$ and $z\ls w$
such that $x+24y$ is a square. We also guess that any $n\in\N$ with $n\not=47$ can be written as $x^2+y^2+z^2+w^2$ $(x,y,z,w\in\N$) with $x+7y$ a square, and
that any integer $n>3$ can be expressed as $x^2+y^2+z^2+w^2$ $(x,y,z,w\in\N$) with $3x-y$ a square.

\proclaim{Conjecture 4.2} {\rm (i)} For each $k=1,2,3$ and $n\in\Z^+$, there are $x,y,w\in\N$ and $z\in\Z^+$ with $n=x^2+y^2+z^2+w^2$ such that $(x+ky)z$ is a square.

{\rm (ii)} Any positive integer can be written $x^2+y^2+z^2+w^2$ with $x,y,w\in\N$ and $z\in\Z^+$ such that
$(ax-by)z$ is a square, where $(a,b)$ is either of the ordered pairs
$$(1,2),\ (2,2),\ (3,2),\ (3,3),\ (4,2),\ (6,6).$$
\endproclaim
\Remark\ 4.2. In view of Theorem 1.4(i), $a(x-y)z\ (a=2,3,6)$, $(x-2y)z$ and $(4x-2y)z$ are all suitable.
We also conjecture that $(ax+by)z$ is suitable for any ordered pair $(a,b)$ among
$$(2,5),\ (3,3),\ (3,6),\ (3,15),\ (5,6),\ (5,11),\ (5,13),\ (5,15),\ (6,15),\ (8,46),\ (9,23).$$

\proclaim{Conjecture 4.3} {\rm (i) (1-3-5-Conjecture)} Any $n\in\N$ can be written as $x^2+y^2+z^2+w^2$ with $x,y,z,w\in\N$
such that $x+3y+5z$ is a square.

{\rm (ii)} Any integer $n > 15$ can be written as $x^2 + y^2 + z^2+w^2$  with $x,y,z,w\in\N$ such that $3x+5y+6z$ is twice a square.

{\rm (iii)} Let $a,b,c\in\Z^+$ with $b\ls c$ and $\gcd(a,b,c)$ squarefree. Then $ax-by-cz$ is suitable if and only if
$(a,b,c)$ is among the five triples
$$(1,1,1),\ (2,1,1),\ (2,1,2),\ (3,1,2),\ (4,1,2).$$

{\rm (iv)} Let $a,b,c\in\Z^+$ with $a\ls b$ and $\gcd(a,b,c)$ squarefree. Then $ax+by-cz$ is suitable if and only if
$(a,b,c)$ is among the following $52$ triples:
$$\align&(1,1,1),\ (1,1,2),\ (1,2,1),\ (1,2,2),\ (1,2,3),\ (1,3,1),
\\&(1,3,3),\ (1,4,4),\ (1,5,1),\ (1,6,6),\ (1,8,6),\ (1,12,4),\ (1,16,1),
\\&(1,17,1),\ (1,18,1),\ (2,2,2),\ (2,2,4),\ (2,3,2),\ (2,3,3),\ (2,4,1),
\\&(2,4,2),\ (2,6,1),\ (2,6,2),\ (2,6,6),\ (2,7,4),\ (2,7,7),\ (2,8,2),
\\&(2,9,2),\ (2,32,2),\ (3,3,3),\ (3,4,2),\ (3,4,3),\ (3,8,3),\ (4,5,4),
\\&(4,8,3),\ (4,9,4),\ (4,14,14),\ (5,8,5),\ (6,8,6),\ (6,10,8),\ (7,9,7),
\\&(7,18,7),\ (7,18,12),\ (8,9,8),\ (8,14,14),\ (8,18,8),\ (14,32,14),
\\&(16,18,16),\ (30,32,30),\ (31,32,31),\ (48,49,48),\ (48,121,48).
\endalign$$
\endproclaim

\Remark\ 4.3. We guess that if $a,b,c$ are positive integers with $\gcd(a,b,c)$ squarefree such that any $n\in\N$ can be written as
 $x^2+y^2+z^2+w^2$ $(x,y,z,w\in\N)$ with $ax+by+cz$ a square then we must have $\{a,b,c\}=\{1,3,5\}$.
Concerning the 1-3-5-Conjecture, see [S16, A271518, A273294, A273302, A278560] for related data; for example,
$43 = 1^2 + 5^2 + 4^2+1^2$ with $1 + 3\times5 + 5\times4 = 6^2$. We have verified parts (i) and (ii) of Conjecture 4.3
for $n$ up to $3\times10^7$. The author would like to offer 1350 US dollars as the prize for the first complete solution of the 1-3-5-Conjecture.
We also conjecture that any $n\in\N$ can be written as $x^2+y^2+z^2+w^2$ with $x,y,z,w\in\Z$ such that $x+3y+5z$ is a
cube. For parts (iii) and (iv) of Conjecture 4.3, see the comments in [S16, A271775].

\proclaim{Conjecture 4.4} {\rm (i)} For each $c=1,2,4$, any $n\in\N$ can be written as $w^2+x^2+y^2+z^2$ with $w,x,y,z\in\N$ and $y\ls z$ such that $2x+y-z=ct^3$ for some $t\in\N$.

{\rm (ii)} Any $n\in\N$ can be written as $w^2+x^2+y^2+z^2$ with $w\in\Z$ and $x,y,z\in\N$ such that $w+x+2y-4z$ is twice a nonnegative cube.

{\rm (iii)} Any $n\in\N$ not of the form $4^{2k+1}\times7\ (k\in\N)$ can be written as $w^2+x^2+y^2+z^2$ $(w,x,y,z\in\N)$ with $w+2x+3y+5z$ a square.
Also, for any $a,b,c,d\in\Z^+$, there are infinitely many positive integers which cannot be written as $w^2 + x^2 + y^2 + z^2$
 $(w,x,y,z\in\N)$ with $aw + bx + cy + dz$ a square.

{\rm (iv)} Let $a,b,c,d\in\Z^+$ with $a\ls b\ls c$ and $\gcd(a,b,c,d)$ squarefree. Then $ax+by+cz-dw$
is suitable if and only if $(a,b,c,d)$ is among the following quadruples
$$\align&(1,1,2,1),\ (1,2,3,1),\ (1,2,3,3),\ (1,2,4,2),\ (1,2,4,4),\ (1,2,5,5),
\\&(1,2,6,2),\ (1,2,8,1),\ (2,2,4,4),\ (2,4,6,4),\ (2,4,6,6),\ (2,4,8,2).
\endalign$$

{\rm (v)} Any positive integer can be written as $w^2+x^2+y^2+z^2$ with $w+x+y-z$ a square, where $w\in\Z$ and $x,y,z\in\N$ with $|w|\ls x\gs y\ls z<x+y$.
Also, any $n\in\N$ can be written as  $w^2+x^2 + y^2 + z^2$ with $w+x+y-z$ a nonnegative cube, where $w,x,y,z$ are integers with $|x|\ls y\gs z\gs0$.

{\rm (vi)} Any $n\in\Z^+$ can be written as $w^2+x^2+y^2+z^2$ with $w\in\Z^+$ and $x,y,z\in\Z$ such that $aw+bx+cy+dz$ is a nonnegative cube, whenever
$(a,b,c,d)$ is among the quadruples
$$(1,1,1,3),\ (1,1,2,2),\ (1,1,2,4),\ (1,1,3,3),\ (8,2,6,8),\ (8,4,4,8),\ (8,4,8,12).$$
\endproclaim
\Remark\ 4.4. See [S16, A273432, A273568, A279522, A272620 and A273458] for related data.

\proclaim{Conjecture 4.5} Let $a,b,c\in\Z^+$ with $a\ls b\ls c$ and $\gcd(a,b,c)$ squarefree.
Then the polynomial $w(ax+by+cz)$ is suitable if and only if $(a,b,c)$ is among the five triples
$$(1,2,3),\ (1,3,6),\ (1,6,9),\ (5,6,9),\ (18,30,114).$$
\endproclaim
\Remark\ 4.5. See [S16, A271724] for the number of ways to write $n=w^2+x^2+y^2+z^2$ with $w,y,z\in\N$ and $x\in\Z^+$ such that $w(x+2y+3z)$ is a square.

\proclaim{Conjecture 4.6} {\rm (i)} Any $n\in\Z^+$ can be written as $w^2+x^2+y^2+z^2$ with $w\in\Z^+$ and $x,y,z\in\N$ such that
$wx+2xy+2yz$
$($or $2wx+xy+4yz)$ is a square. Also, any $n\in\Z^+$ can be written as $w^2+x^2+y^2+z^2$ with $w\in\Z^+$, $x,y,z\in\N$ and $x\ls y$ such that
$2xy+yz-zw-wx$ is a square.

{\rm (ii)} Any $n\in\Z^+$ can be written as $w^2+x^2+y^2+z^2$ with $w\in\Z^+$ and $x,y,z\in\N$ such that
$w^2+4xy+8yz+32zx$ is a square.

{\rm (iii)} For each $k = 1, 2, 8,16,48,336$, any positive integer can be written as $w^2 + x^2 + y^2 + z^2$
 with $w\in\Z^+$ and $x,y,z\in\N$ such that $w^2 + k(xy+yz)$ is a square.

{\rm (iv)} Let $a,b,c\in\Z^+$ with $\gcd(a,b,c)$ squareferee. Then the polynomial $axy+byz+czx$ is suitable if and only if $\{a,b,c\}$ is among
$$\{1,2,3\},\ \{1,3,8\},\ \{1,8,13\},\ \{2,4,45\},\ \{4,5,7\},\ \{4,7,23\},\ \{5,8,9\},\ \{11,16,31\}.$$

{\rm (v)} Any positive integer can be written as $w^2 + x^2 + y^2 + z^2$ with $w\in\Z^+$ and $x,y,z\in\N$ such that
 $wx + xy + 2yz + 3zx$
 $($or $wx + 3xy + 8yz + 5zx)$ is twice a square. Also, each $n\in\Z^+$ can be written as $w^2 + x^2 + y^2 + z^2$ with $w\in\Z^+$ and $x,y,z\in\N$ such that
 $6wx + 2xy + 3yz + 4zx=3t^2$ for some $t\in\N$.
\endproclaim
\Remark\ 4.6. See [S16, A271644, A273021 and A271665] for related data.

\proclaim{Conjecture 4.7} {\rm (i)} Any natural number can be written as $x^2 + y^2 + z^2 + w^2$ with $x,y,z,w\in\N$ and $x\gs y$ such that
$ax^2+by^2+cz^2$ is a square, provided that the triple $(a,b,c)$ is among
$$(1,8,16),\ (4,21,24),\ (5,40,4),\ (9,63,7),\ (16,80,25),\ (36,45,40),\ (40,72,9).$$

{\rm (ii)} $ax^2+by^2+cz^2$ is suitable if $(a,b,c)$ is among the triples
$$\gather (1,3,12),\ (1,3,18),\ (1,3,21),\ (1,3,60),\ (1,5,15),\ (1,8,24),\ (1,12,15),\ (1,24,56),
\\(3,4,9),\ (3,9,13),\ (4,5,12),\ (4,5,60),\ (4,9,60),\ (4,12,21),\ (4,12,45),\ (5,36,40).
\endgather$$

{\rm (iii)} If $a,b,c$ are positive integers with $ax^2+by^2+cz^2$ suitable, then $a,b,c$  cannot be pairwise coprime.
\endproclaim
\Remark\ 4.7. See [S16, A271510 and A271513] for related data.

\proclaim{Conjecture 4.8} {\rm (i)} Any $n\in\Z^+$ can be written as $w^2 + x^2 + y^2 + z^2$ with $w\in\Z^+$ and $x,y,z\in\N$
such that $(10w+5x)^2+(12y+36z)^2$ is a square.

{\rm (ii)} Each positive integer can be written as $x^2+y^2+z^2+w^2$ with $x,y,z,w\in\N$ and $y>z$ such that $(x+y)^2+(4z)^2$ is a square.

{\rm (iii)} Any integer $n>5$ can be written as $x^2+y^2+z^2+w^2$ with $x,y,z,w\in\N$, $x+y>0$ and $z>0$  such that
$(8x+12y)^2+(15z)^2$ is a square (i.e., $8x+12y$ and $15z$ are the two legs of a right triangle
with positive integer sides).

{\rm (iv)} Any $n\in\N$ can be written as $x^2+y^2+z^2+w^2$ with $x,y,z,w\in\N$ and $x+y\gs z$ such that
$(x+y+z)^2+(4(x+y-z))^2$ is a square.

{\rm (v)} Any $n\in\Z^+$ can be written as $x^2+y^2+z^2+w^2$ with $x,y,z,w\in\N$ and $y<z$ such that
$x+8y+8z+15w$ and $6(x+y+z+w)$ are the two legs of a right triangle
with positive integer sides. Also, all the polynomials
$$\align &(x+3y+6z+17w)^2+(20x+4y+8z+4w)^2,
\\&(x+3y+9z+17w)^2+(20x+4y+12z+4w)^2,
\\&(x+3y+11z+15w)^2+(12x+4y+4z+20w)^2,
\\&(3(x+2y+3z+4w))^2+(4(x+4y+3z+2w))^2,
\\&(3(x+2y+3z+4w))^2+(4(x+5y+3z+w))^2
\endalign$$
are suitable.
\endproclaim
\Remark\ 4.8. This conjecture is particularly mysterious since it is related to Pythagorean triples. See [S16, A271714, A273108, A273107 and A273134] for related data.

\proclaim{Conjecture 4.9} {\rm (i)} Any $n\in\Z^+$ can be written as $w^2+x^2(1+y^2+z^2)$ with $w,x,y,z\in\N$, $x>0$ and $y\eq z\pmod 2$.
Moreover, any $n\in\Z^+$ with $n\not=449$ can be written as $4^k(1+x^2+y^2)+z^2$ with $k,x,y,z\in\N$ and $x\eq y\pmod 2$.

{\rm (ii)} Each $n\in\Z^+$ can be written as $4^k(1+x^2+y^2)+z^2$ with $k,x,y,z\in\N$ and $x\ls y\ls z$.

{\rm (iii)} Any $n\in\Z^+$ can be written as $4^k(1+5x^2+y^2)+z^2$ with $k,x,y,z\in\N$, and also
each $n\in\Z^+$ can be written as $4^k(1+x^2+y^2)+5z^2$ with $k,x,y,z\in\N$.
\endproclaim
\Remark\ 4.9. This conjecture was motivated by Theorem 1.2(v). We can show the first part provided that any integer $n>432$ can be written as $2x^2+2y^2+z(4z+1)$ with $x,y,z\in\Z$.
Under the GRH we are able to prove Conjecture 4.9(iii) with the help of the work in [KS]. See [S16, A275738, A275656, A275675 and A275676] for related data.

\proclaim{Conjecture 4.10} {\rm (i)} $x^2+kyz$ $(k=12,24,32,48,84,120,252)$, $9x^2-4yz$, $9x^2+140yz$, $25x^2+24yz$ and $121x^2-20yz$ are all suitable.

{\rm (ii)} The polynomials $w(x^2+8y^2-z^2)$, $(3x^2+13y^2)z$, $(5x^2+11y^2)z$, $(15x^2+57y^2)z$, $(15x^2+165y^2)z$ and $(138x^2+150y^2)z$
are all suitable.

{\rm (iii)} Any positive integer can be written as $w^2+x^2+y^2+z^2$ with $w\in\Z^+$ and $x,y,z\in\N$ such that $36x^2y+12y^2z+z^2x$
$($or $x^3+4yz(y-z)$, or $x^3+8yz(2y-z))$ is a square.

{\rm (iv)} Let $a$ and $b$ be nonzero integers with $\gcd(a,b)$ squarefree. Then the polynomial $ax^4+by^3z$ is suitable if and only if $(a,b)$ is among the ordered pairs
$$(1,1),\ (1,15),\ (1,20),\ (1,36),\ (1,60),\ (1,1680),\ (9,260).$$

{\rm (v)} For each triple $(a,b,c)=(1,20,60),\ (1,24,56),\ (9,20,60),\ (9,32,96)$, any $n\in\Z^+$ can be written as $x^2+y^2+z^2+w^2$
with $x,y,z\in\N$ and $w\in\Z^+$ such that $ax^4+by^3z+cyz^3$ is a square.
\endproclaim
\Remark\ 4.10. See [S16, A272888, A272332, A279056, A272336, A280831 and A272351] for related data and comments. Part (iii) of Conjecture 4.10 looks very curious.
The author ever guessed that $x^2-4yz$, $x^2+4yz$ and $x^2+8yz$ are suitable, then his student You-Yin Deng made a clever observation:
$$x^2+2y^2+6z^2=x^2+(y+z)^2+(y-z)^2+(2z)^2\ \t{with}\ (2z)^2+4(y+z)(y-z)=(2y)^2.$$
Since $x^2\pm4yz=x^2\pm(y+z)^2\mp(z-y)^2$, and $x(x+y-z)$ is suitable (cf. Theorem 1.4(i)), $x^2+4yz$ and $x^2-4yz$ are indeed suitable.
The author's student Yu-Chen Sun has proved that $x^2+8yz$ is also suitable.

\proclaim{Conjecture 4.11} {\rm (i)} Any positive integer can be written as $x^2+y^2+z^2+w^2$ with $x,y,z,w\in\N$ and $z<w$ such that
$4x^2+5y^2+20zw$ is a square. Also, the polynomials $x^2+8y^2+8zw$, $(3x+5y)^2-24zw$, $(x-2y)^2+24zw$ and $(x-3y)^2+16zw$ are suitable.

{\rm (ii)} The polynomials $x^2+3y^2+4z^2+(x+y+z)^2$, $x^2 + 3y^2 + 5z^2 - 8w^2$, $(x-2y)^2+8z^2+16w^2$, $4(x-3y)^2+9z^2+12w^2$, $(x+2y)^2+8z^2+40w^2$ and $9(x+2y)^2+16z^2+24w^2$ are all suitable.

{\rm (iii)} The polynomial $w^2x^2+3x^2y^2+2y^2z^2$ is suitable.

{\rm (iv)} Any positive integer can be written as $w^2+x^2+y^2+z^2$ with $w\in\Z^+$ and $x,y,z\in\N$ such that
$w^2x^2+5x^2y^2+80y^2z^2+20z^2w^2$ is a square.
\endproclaim
\Remark\ 4.11. See [S16, A272084, A271778, A271824, A273278, A269400 and A262357] for related data and more similar conjectures.

\proclaim{Conjecture 4.12} {\rm (i)} Let $a,b,c,d\in\Z^+$ with $a\ls b$ and $c\ls d$, and $\gcd(a,b,c,d)$ squarefree.
 Then $ax+by-cz-dw$ is suitable if and only if $(a,b,c,d)$ is among the quadruples
 $$\gather (1,2,1,1),\ (1,2,1,2),\ (1,3,1,2),\ (1,4,1,3),
 \\ (2,4,1,2),\ (2,4,2,4),\ (8,16,7,8),\ (9,11,2,9),\ (9,16,2,7).
 \endgather$$

{\rm (ii)} The polynomial
$ax^2+by^2-cz^2-dw^2$ is suitable if $(a,b,c,d)$ is among the quadruples
$$\align&(3,9,3,20),\ (5,9,5,20),\ (5,25,4,5),\ (9,81,9,20),\ (12,16,3,12),\ (16,64,15,16),
\\&(20,25,4,20),\ (27,81,20,27),\ (30,64,15,30),\ (32,64,15,32),\ (48,64,15,48).
\endalign$$
\endproclaim
\Remark\ 4.12. We also conjecture that $ax^2-by^2-cz^2$ is suitable if $(a,b,c)$ is among the triples
 $$(21,5,15),\ (36,3,8),\ (48,8,39),\ (64,7,8),\ (40,15,144),\ (45,20,144),\ (69,20,60).$$

\proclaim{Conjecture 4.13} {\rm (i)} The polynomial $xyz(x+9y+11z+10w)$ is suitable.

{\rm (ii)} Any $n\in\Z^+$ can be written as $x^2+y^2+z^2+w^2$ with $x\in\Z^+$, $y,z,w\in\N$ and $y\gs z$
such that $xyz(x+3y+13z)$ is a square.

{\rm (iii)} Any $n\in\Z^+$ can be written as
$x^2+y^2+z^2+w^2$ with $xy(3x+5y+2z+3w)$
$($or $xy(x+11y+z+2w)$, or $2xy(x+2y+z+2w)$, or $2xy(x+6y+z+2w))$ a square,
where $x,y,z,w$ are nonnegative integers with $w>0$ $($or $z>0)$.

{\rm (iv)} The polynomial $xy(ax^2+by^2+cz^2)$ is suitable
whenever $(a,b,c)$ is among the triples
$$\align&(1,8,20),\ (3,5,15),\ (6,14,4),\ (7,9,5),\ (7,29,5),\ (18,38,18),\ (39,81,51).
\endalign$$

{\rm (v)} If $(a,b,c)$ is one of the six triples
$$(1,2,4),\ (1,2,9),\ (1,3,4),\ (2,3,4),\ (2,4,6),\ (4,8,10),$$
then any $n\in\Z^+$ can be written as $w^2+x^2+y^2+z^2$ with $w\in\Z^+$ and $x,y,z\in\N$
such that $w(25w+24(ax+by+cz))$ is a square.
\endproclaim
\Remark\ 4.13. See [S16,  A267121, A260625, A261876 and A268197] for related data.

\proclaim{Conjecture 4.14} {\rm (i)} Any $n\in\N$ can be written as $x^2+y^2+z^2+w^2$ with $x,y,z,w\in\N$ such that $xy+2zw$ or $xy-2zw$ is a square.
Also, each $n\in\N$ can be written as $x^2+y^2+z^2+w^2$ with $x,y,z,w\in\N$ and $\max\{x,y\}\gs\min\{z,w\}$
such that $xy+zw/2$ or $xy-zw/2$ is a square.

{\rm (ii)} Any $n\in\N$ can be written as $x^2+y^2+z^2+w^2$ with $x,y,z\in\N$, $w\in\Z$ and $x\gs z$ such that $3x^2y+z^2w$ is a square. Also, for each ordered pair $(a,b)=(7,1),\ (8,1),\ (9,2)$,
any $n\in\N$ can be written as $x^2+y^2+z^2+w^2$ with $x,y,z\in\N$ and $w\in\Z$
such that $ax^2y+bz^2w$ is a square.

{\rm (iii)}  Any $n\in\Z^+$ can be written as $x^2+y^2+z^2+w^2$ with $x\in\Z^+$, $y\in\N$ and $z,w\in\Z$
such that $xy+yz+zw$ is a fourth power. Also,
any $n\in\Z^+$ can be written as $x^2+y^2+z^2+w^2$ with $x,y,z\in\Z$ and $w\in\Z^+$
such that $xy+yz+2zw+2wx$ is a fifth power.
\endproclaim
\Remark\ 4.14. See [S16, A270073, A272977 and A273826] for related data. We have verified Conjecture 4.14(i) for all $n=0,1,\ldots,2\times10^5$.

\proclaim{Conjecture 4.15} Any $n\in\N$ can be written as $p_5(u)+p_5(v)+p_5(x)+p_5(y)+p_5(z)$ with $u,v,x,y,z\in\N$ such that
$u+2v+4x+5y+6z$ is a pentagonal number, where $p_5(k)$ with $k\in\N$ denotes the pentagonal number $k(3k-1)/2$.
\endproclaim
\Remark\ 4.15. As conjectured by Fermat and proved by Cauchy, each natural number can be written as the sum of five pentagonal numbers
(cf. [N96, pp.\,27-34] or [MW. pp.\,54-57]). See [S16, A271608] for some data related to Conjecture 4.15.

\proclaim{Conjecture 4.16} {\rm (i)} Any $n\in\N$ can be written as $\sum_{i=1}^9x_i^3$ with $x_i\in\N$
such that
$$x_1+x_2+x_3+2x_4+3x_5+4x_6+4x_7+9x_8+15x_9$$
is a cube.

{\rm (ii)} Any $n\in\N$ can be written as $\sum_{i=1}^9x_i^3$ with $x_i\in\N$
such that
$$x_1^3+x_2^3+x_3^3+2x_4^3+3x_5^3+4x_6^3+5x_7^3+14x_8^3+19x_9^3$$
is a cube.

{\rm (iii)} Any $n\in\N$ can be written as $\sum_{i=1}^9x_i^3$ with $x_i\in\N$
such that $\sum_{i=1}^9ix_i^4$
$($or $\sum_{i=1}^9ix_i^2)$ is a square.
\endproclaim
\Remark\ 4.16. It is well-known that any natural number is the sum of nine nonnegative cubes (cf. [N96, pp.\,41-43]).
We even conjecture further that any $n\in\N$ can be written as $u^3+v^3+2x^3+2y^3+3z^3$ with $u,v,x,y,z\in\N$ (cf. [S16, A271099]).

\medskip
\Ack. The author would like to thank his students You-Yin Deng and Yu-Chen Sun, and the annonymous referee for helpful comments.

\enddocument